\documentclass{amsart}
\usepackage[a4paper]{geometry}
\usepackage{amsmath}
\usepackage{amssymb}
\usepackage{amsthm}
\usepackage{mathrsfs}
\usepackage{mathtools}

\usepackage{tikz}
\usepackage{cleveref}
\usepackage{tikz-cd}
\usepackage{mathtools}

\usepackage{placeins}

\usepackage{float}

\usepackage{comment}

\usepackage{enumerate}

\usepackage{xinttools} 
\usetikzlibrary{calc}

\theoremstyle{definition}
\newtheorem{thm}{Theorem}[section]
\crefname{thm}{Theorem}{Theorems}
\newtheorem{cor}[thm]{Corollary}
\newtheorem{prop}[thm]{Proposition}
\crefname{prop}{Proposition}{Propositions}

\crefname{lem}{Lemma}{Lemmas}

\crefname{defn}{Definition}{Definitions}

\newtheorem{obs}[thm]{Observation}

\newtheorem*{ack*}{Acknowledgements}

\newcommand{\card}{\#}

\title{Sharp bounds for the Tao-Vu Discrete John's Theorem}
\author{Peter van Hintum}
\thanks{New College, University of Oxford, UK. email: peter.vanhintum@new.ox.ac.uk}
\author{Peter Keevash}
\thanks{Mathematical Institute, University of Oxford, UK. Supported by ERC Advanced Grant 883810.}
\allowdisplaybreaks

\begin{document}
\maketitle
\begin{abstract}
Tao and Vu showed that every centrally symmetric convex progression $C\subset\mathbb{Z}^d$ 
is contained in a generalised arithmetic progression of size $d^{O(d^2)} \card C$. 
Berg and Henk improved the size bound to $d^{O(d\log d)} \card C$. 
We obtain the bound $d^{O(d)} \card C$, which is sharp up to the implied constant,
and is of the same form as the bound in the continuous setting given by John's Theorem.
\end{abstract}

\section{Introduction}

A classical theorem of John \cite{fritz1948extremum} shows that for any centrally symmetric convex set $K\subset \mathbb{R}^d$
there exists an ellipsoid $E$ centred at the origin so that $E\subset K\subset \sqrt{d}E$. 
This immediately implies that there exists a parallelotope $P$ so that $P\subset E\subset K\subset \sqrt{d}E\subset dP$. 
In the discrete setting, quantitative covering results are of great interest in Additive Combinatorics,
a prominent example being the Polynomial Freiman-Ruzsa Conjecture, which asks for effective bounds
on covering sets of small doubling by convex progressions. 
In this context, a natural analogue of John's Theorem in $\mathbb{Z}^d$ would be
covering centrally symmetric convex progressions by generalised arithmetic progressions. 
Here, a $d$-dimensional \emph{convex progression} is a set of the form $K\cap \mathbb{Z}^d$, 
where $K\subset \mathbb{R}^d$ is convex, and a $d$-dimensional \emph{generalized arithmetic progression} ($d$-GAP) 
is a translate of a set of the form $\left\{\sum_{i=1}^d m_ia_i: 1\leq m_i\leq n_i\right\}$ for some $n_i\in\mathbb{N}$ and $a_i\in \mathbb{Z}^d$. 

Tao and Vu \cite{tao2006additive,tao2008john} obtained such a discrete version of John's Theorem,
showing that for any centrally symmetric  $d$-dimensional convex progression $C\subset\mathbb{Z}^d$
there exists a $d$-GAP $P$ so that $P\subset C\subset O(d)^{3d/2}\cdot P$, 
where $m\cdot P:=\left\{\sum_{i=1}^m p_i: p_i\in P\right\}$ denotes the iterated sumset.
Berg and Henk \cite{berg2019discrete} improved this to $P\subset C\subset d^{O(\log(d))}\cdot P$.
Our focus will be on the covering aspect of these results, i.e.\ minimising the ratio $\card P' / \card C$,
where $P'$ is a $d$-GAP covering $C$. This ratio is bounded by $d^{O(d^2)}$ by Tao and Vu
and by $d^{O(d\log d)}$ by Berg and Henk. We obtain the bound $d^{O(d)}$, which is optimal.

\begin{thm}\label{mainthm}
For any centrally symmetric  convex progression $C\subset\mathbb{Z}^d$
there exists a $d$-GAP $P$ containing $C$ with $\card P\leq O(d)^{3d} \card C$.
\end{thm}

\begin{cor}\label{maincor}
For any centrally symmetric  convex progression $C\subset\mathbb{Z}^d$ and linear map $\phi:\mathbb{R}^d\to\mathbb{R}$,
there exists a $d$-GAP $P$ containing $C$ with $\card \phi(P)\leq O(d)^{3d} \card \phi(C)$.
\end{cor}

The optimality of \Cref{mainthm} is demonstrated by the intersection of a ball $B$ with a lattice $L$. 
Moreover, Lovett and Regev \cite{lovett2017} obtained a more emphatic negative result,
disproving the GAP analogue of the Polynomial Freiman-Ruzsa Conjecture,
by showing that by considering a random lattice $L$
one can find a convex $d$-progression  $C = B \cap L$ such that any $O(d)$-GAP $P$
with $|P| \le |C|$ has $|P \cap C| < d^{-\Omega(d)} |C|$. Our result can be viewed 
as the positive counterpart that settles this line of enquiry, showing that indeed $d^{\Theta(d)}$
is the optimal ratio for covering convex progressions by GAPs.
 
\section{Proof}

We start by recording two simple observations and a property of Mahler Lattice Bases.

\begin{obs}\label{ContJohn}
Given a centrally symmetric convex set $K\subset\mathbb{R}^d$, there exists a centrally symmetric parallelotope $Q$ and a centrally symmetric ellipsoid $E$ so that $\frac1d Q\subset E\subset  K\subset\sqrt{d}E\subset Q$, so in particular $|Q|\leq d^{d}|K|$.
\end{obs}

\begin{obs}\label{anybasistransformation}
Let $X,X'\in\mathbb{R}^{d\times d}$ be so that the rows of $X$ and $X'$ generate the same lattice of full rank in $\mathbb{R}^d$.
Then $\exists T\in GL_n(\mathbb{Z})$ so that $TX=X'$.
\end{obs}

\begin{prop}[Corollary 3.35 from \cite{tao2006additive}]\label{taovubasisprop}
Given a lattice $\Lambda\subset\mathbb{R}^d$ of full rank, there exists a lattice basis $v_1,\dots, v_d$ of $\Lambda$ 
so that $\prod_{i=1}^{d} \|v_i\|_2 \leq O(d^{3d/2})\det(v_1,\dots, v_d)$.
\end{prop}

With these three results in mind, we prove the theorem.

\begin{proof}[Proof of \Cref{mainthm}]
By passing to a subspace if necessary, we may assume that $C$ is full-dimensional.
Write $C = K \cap \mathbb{Z}^d$ where $K\subset\mathbb{R}^d$ is centrally symmetric and convex. 
Use \Cref{ContJohn} to find a parallelotope $Q\supset K$ so that $|Q|\leq d^d |K|$. 
Let the defining vectors of $Q$ be $u_1,\dots,u_d$, i.e. $Q=\{\sum \lambda_i u_i: \lambda_i\in[-1,1]\}$. 
Write $u_i^j$ for the $j$-th coordinate of $u_i$ and write $U$ for the matrix $(u_i^j)$ with rows $u^j$ and columns $u_i$. 

Consider the lattice $\Lambda$ generated by the vectors $u^j$ (these are the vectors formed by the $j$-th coordinates of the vectors $u_i$). 
Using \Cref{taovubasisprop} find a basis $v^1,\dots,v^d$ of $\Lambda$ so that $\prod_{j=1}^{d} |v^j|\leq O(d^{3d/2})\det(v^1,\dots, v^d)$. 
Write $v^j_i$ for the $i$-th coordinate of $v^j$ and write $V:=(v_i^j)$. By \Cref{anybasistransformation}, we can find $T\in GL_n(\mathbb{Z})$ so that $TU=V$. 
Let $T'\colon\mathbb{R}^d\mapsto\mathbb{R}^d$ defined by $T' u_i = v_i$ for $1 \le i \le d$.
Note that $T'(\mathbb{Z}^d)=\mathbb{Z}^d$ and in the standard basis $T'$ corresponds to matrix multiplication by $T$.

Write $Q':=T'(Q)=\{\sum \lambda_i v_i: \lambda_i\in[-1,1]\}$ and consider the smallest axis aligned box $B:=\prod [-a_i,a_i]$ containing $Q'$. 
Note that $a_j\leq\sum_{i} |v_i^j|=||v^{j}||_1\leq \sqrt{d}||v^{j}||_2$. Hence, we find
\begin{align*}
|B|&=2^d\prod_{i=1}^d a_i\leq 2^d\prod_{j=1}^d \sqrt{d}||v^j||_2\leq O(d)^{2d} \det(v^1,\dots,v^d)= O(d)^{2d} \det(v_1,\dots,v_d)=O(d)^{2d} |Q'|.
\end{align*}
Now we cover $C$ by a $d$-GAP $P$, constructed by the following sequence:
$$C=K\cap \mathbb{Z}^d\subset Q\cap \mathbb{Z}^d= T'^{-1}(Q')\cap \mathbb{Z}^d
\subset T'^{-1}(B)\cap \mathbb{Z}^d = T'^{-1}(B\cap \mathbb{Z}^d) =: P.$$
It remains to bound $\card P$. As $C$ is full-dimensional each $a_i \ge 1$, so
\begin{align*}
\card P&=\card (B \cap  \mathbb{Z}^d) \leq 2^d |B|\leq O(d)^{2d} |Q'|= O(d)^{2d} |Q|
\leq O(d)^{3d} |K|\leq O(d)^{3d} \card C. \qedhere
\end{align*}
\end{proof}

\begin{proof}[Proof of \Cref{maincor}]
Let $m:=\max_{x\in\mathbb{Z}}\card(\phi^{-1}(x)\cap C)$ and note that
$\card \phi(C)\geq \card C/ m$. Analogously, let $m':=\max_{x\in\mathbb{Z}}\card(\phi^{-1}(x)\cap P)$, so that $m'\geq m$. By translation, we may assume that $m'$ is achieved at $x=0$.
Note that for any $x = \phi(p)$ with $p \in P$ and $p' \in P \cap \phi^{-1}(0)$ we have $p+p' \in P+P$ with $\phi(p+p')=x$, so $\card(\phi^{-1}(x)\cap (P+P))\geq m'$. We conclude that
\begin{align*}
\card \phi(P) &\leq \card (P+P)/m'\leq 2^d \card P/m\leq O(d)^{3d}\card C/m\leq O(d)^{3d}\card \phi(C). \qedhere
\end{align*}
\end{proof}

\bibliographystyle{alpha}
\bibliography{references}

\end{document}